\documentstyle[11pt,amstex,amscd]{article}

\newtheorem{theorem}{Theorem}[section]
\newtheorem{proposition}{Proposition}[section]
\newtheorem{lemma}{Lemma}[section]
\newtheorem{corollary}{Corollary}[section]
\newtheorem{definition}{\it Definition}[section]

\newtheorem{remark}{\it Remark}[section]
\newtheorem{notation}{\it Notation}[section]
\newtheorem{claim}{Claim}[section]
\newtheorem{question}{Question}[section]

\newcommand{\qed}{\hspace*{\fill} 
\setlength{\unitlength}{1mm}
\begin{picture}(2.5,2.5)
      \put(0,0){\framebox(2.5,2.5){}}
  \end{picture}
\setlength{\unitlength}{1pt}}

\newcommand{\bs}{\backslash}
\newcommand{\GmG}{{G/\Gamma}}

\newcommand{\proof}{\noindent{\it Proof.~}}

\newcommand{\abscont}{\ll} 

\newcommand{\cl}[1]{\overline{#1}}

\newcommand{\la}[1]{{\frak{ \lowercase{#1}}}}

\newcommand{\inv}{{^{-1}}}
 
\newcommand{\Aut}{{\rm Aut}}
\newcommand{\Aff}{{\rm Aff}}

\newcommand{\End}{{\rm End}}

\newcommand{\rank}{{\rm rank}}

\newcommand{\Ad}{{\rm Ad}}

\def\supp{{\rm supp}}

\def\pr{{\rm pr}}

\def\NN{{\Bbb N}}
\def\RR{{\Bbb R}}
\def\ZZ{{\Bbb Z}}
\def\QQ{{\Bbb Q}}

\def\cB{{\cal B}}

\def\cH{{\cal H}}
\def\cI{{\cal I}}
\def\cJ{{\cal J}}

\def\cP{{\cal P}}
\def\cQ{{\cal Q}}

\def\bfG{{\bf G}}

\def\bfa{{\bf a}}
\def\bfb{{\bf b}}
\def\bfc{{\bf c}}

\def\bfe{{\bf e}}

\def\bfp{{\bf p}}

\def\bft{{\bf t}}

\def\bfv{{\bf v}}
\def\bfw{{\bf w}}

\title{Limit distributions of expanding translates of certain  
orbits on homogeneous spaces}

\author{Nimish A. Shah}

\date{}

\begin{document}


\noindent{\bf\large Limit distributions of expanding  translates of
certain orbits on \\[4pt] homogeneous spaces}

\bigskip

\oddsidemargin  2.0cm   
\textwidth 14.0cm    

\noindent NIMISH A SHAH\\[3pt]
{\small School of Mathematics, Tata Institute of Fundamental
Research, Homi Bhabha Road,\\ Bombay 400 005. E-mail:
nimish@@tifrvax.tifr.res.in}

\bigskip\noindent
{\small MS Received 21 May 1995; revised 9 November 1995}

\bigskip\noindent {\small {\bf Abstract.}  Let $L$ be a Lie group and
$\Lambda$ a lattice in $L$. Suppose $G$ is a non-compact simple Lie
group realized as a Lie subgroup of $L$ and
$\overline{G\Lambda}=L$. Let $a\in G$ be such that $\mbox{Ad}a$ is
semisimple and not contained in a compact subgroup of
$\mbox{Aut}(\mbox{Lie}(G))$.  Consider the expanding horospherical
subgroup of $G$ associated to $a$ defined as $U^+=\{g\in G:
a^{-n}ga^n\to e\mbox{ as } n\to\infty\}$. Let $\Omega$ be a nonempty
open subset of $U^+$ and $n_i\to\infty$ be any sequence. It is showed
that $\overline{\cup_{i=1}^\infty a^{n_i}\Omega\Lambda}=L$. A stronger
measure theoretic formulation of this result is also obtained. Among
other applications of the above result, we describe $G$-equivariant
topological factors of $L/\Lambda\times G/P$, where the real rank of
$G$ is greater than $1$, $P$ is a parabolic subgroup of $G$ and $G$
acts diagonally. We also describe equivariant topological factors of
unipotent flows on finite volume homogeneous spaces of Lie groups.}


\oddsidemargin  0cm   
\textwidth 16.0cm    

\section{Introduction}
Let $G$ be a connected semisimple Lie group with no compact factors
and of $\RR$-rank $\geq 2$, $P$ a parabolic subgroup of $G$, and
$\Gamma$ an irreducible lattice in $G$. It was proved by
Margulis~\cite{Mar:factors} that if $\phi:G/P\to Y$ is a measure class
preserving $\Gamma$-equivariant factor of $G/P$ then there exist a
parabolic subgroup $Q$ containing $P$ and a measurable isomorphism
$\psi:Y\to G/Q$ such that $\psi\circ\phi$ is the canonical quotient
map. The topological analogue of this result was obtained by
Dani~\cite{Dani:factors}, who proved that, in the above notation, if
$\phi$ is continuous then $\psi$ can be chosen to be a
homeomorphism. On the other hand in the result of Margulis was
generalized by Zimmer~\cite{Zim:factors} in the measure theoretic
category. This result was later used in \cite{Stuck:Zimmer} for
describing faithful and properly ergodic finite measure preserving
$G$-actions. It was suggested by Stuck~\cite{Stuck:msri} that the
following question, which is a topological analogue of Zimmer's
result, is of importance for studying locally free minimal
$G$-actions. 

\begin{question}   \label{que:stuck}
Let $G$ be a simple Lie group of $\RR$-rank$\geq 2$. Suppose that $G$
acts minimally and locally freely on a compact Hausdorff space
$X$. Suppose there are $G$-equivariant continuous surjective maps
$X\times G/P\stackrel{\phi}{\to} Y \stackrel{\psi}{\to} X$ such that
$\psi\circ\phi$ is the projection onto $X$, where $G$ acts diagonally
on $X\times G/P$. Does there exist a parabolic subgroup $Q$ containing
$P$ and a $G$-equivariant homeomorphism $\rho:Y\to X\times G/Q$ such
that $\rho\circ\phi$ is the canonical quotient map?
\end{question}

The above mentioned result of Dani says that this question has the
affirmative answer if $X=G/\Gamma$, $\Gamma$ being a lattice in
$G$. In this paper we consider the case when $G$ is a Lie subgroup of
a Lie group $L$ acting on $X=L/\Lambda$ by translations, $\Lambda$
being a lattice in $L$. To analyze this case we follow the method of
the proof of Dani~\cite{Dani:factors}. To adapt Dani's proof for the
general case one needs the following theorem~\ref{thm:parabolic},
which is a nontrivial generalization of its particular case of $L=G$
(cf.~\cite[Lemma~1.1]{Dani:factors}). Its proof involves, in an
essential way, Ratner's theorem~\cite{R:measure} on classification of
finite ergodic invariant measures of unipotent flows on homogeneous
spaces.

For the results stated in the introduction, let $L$ denote a connected
Lie group, $\Lambda$ a lattice in $L$, $\pi:L\to L/\Lambda$ the
natural quotient map, and $\mu_L$ the (unique) $L$-invariant
probability measure on $L/\Lambda$. 

\begin{theorem}  \label{thm:parabolic}
Let $G$ be a connected semisimple Lie group. Let $a\in G$ be a
semi-simple element; that is, $\Ad(a)$ is a semi-simple endomorphism of
the Lie algebra of $G$. Consider the expanding horospherical subgroup
$U^+$ of $G$ associated to $a$ which is defined as 
\[
U^+=\{u\in G: a^{-n} u a^n \to e \mbox{ as $n\to\infty$}\}.
\]
Assume that $U^+$ is not contained in any proper closed normal subgroup of
$G$. 

Suppose that $G$ is realized as a Lie subgroup of $L$ and
that $\cl{\pi(G)}=L/\Lambda$. Then 
\[
\cl{\pi\left(\cup_{n=1}^\infty a^n U^+\right)}=L/\Lambda.
\]

In particular, if $P$ is any parabolic subgroup of $G$ and
$\cl{\pi(G)}=L/\Lambda$,  then
$\cl{\pi(P)}=L/\Lambda$. 
\end{theorem}

In the case of $L=G$ this result is well-known
(see~\cite[Prop.~1.5]{DR:frames}). Actually
theorem~\ref{thm:parabolic} is a straightforward consequence of a
technically much stronger result stated later in the introduction as
theorem~\ref{thm:translate:horosphere}.

Using the techniques of \cite{Dani:factors} along with
theorem~\ref{thm:parabolic} and the result of
Ratner~\cite{R:equidistribution} on closures of orbits of unipotent flows
on finite volume homogeneous spaces, in the 
next result we provide an affirmative answer to Question~\ref{que:stuck}
in case when $X=L/\Lambda$. In this case we are able relax certain other
conditions in the question as well. 

\begin{theorem}  \label{thm:G/P}
Let $G$ be a semisimple Lie group of $\RR$-$\rank\geq 2$ and with
finite center. Suppose that $G$ is realized as a Lie subgroup of $L$
such that the $G$-action is ergodic with respect to $\mu_L$, and that
$\cl{G_1x}=\cl{Gx}$ for any $x\in L/\Lambda$ and any closed normal
connected subgroup $G_1$ of $G$ such that
$\RR$-$\rank(G/G_1)\leq1$. Let $P$ be a parabolic subgroup of $G$ and
consider the diagonal action of $G$ on $L/\Lambda \times G/P$.  Let
$Y$ be a Hausdorff space with a continuous $G$-action and
$\phi:L/\Lambda \times G/P \to Y$ a continuous $G$-equivariant map.
Then there exist a parabolic subgroup $Q\supset P$, a locally compact
Hausdorff space $X$ with a continuous $G$-action, a continuous
surjective $G$-equivariant map $\phi_1:L/\Lambda \to X$, and a
continuous $G$-equivariant map $\psi:X\times G/Q\to Y$ such that the
following holds:
\begin{enumerate}

\item 
If we define  $\rho:L/\Lambda \times G/P \to X\times G/Q$ as
$\rho(x,gP)=(\phi_1(x),gQ)$ for all $x\in L/\Lambda$ and 
$g\in G$, then  
\[
\phi=\psi\circ\rho.
\]

\item 
There exists an open dense $G$-invariant set $X_0\subset L/\Lambda$ such
that if we put $Z_0=\phi_1(X_0)\times G/Q$ and $Y_0=\psi(Z_0)$, then
$Z_0=\psi\inv(Y_0)$ and $\psi|_{Z_0}$ is injective.  

Furthermore if $Y$ is a locally compact second countable space and
$\phi$ is surjective, then $Y_0$ is open and dense in $Y$ and
$\psi|_{Z_0}$ is a homeomorphism onto $Y_0$.
\end{enumerate}
\end{theorem}

In the next result we classify the $G$-equivariant factors of $L/\Lambda$,
in particular we describe the factor $\phi_1:L/\Lambda\to X$ appearing in
the statement of theorem~\ref{thm:G/P}. The proof of this result uses
the theorem of Ratner on orbit closures of unipotent flows and the main
result of \cite{MS:limit}. 

\begin{definition} \label{def:affine} \rm
Let $\Lambda_1$ be a closed subgroup of $L$. A homeomorphism $\tau$ on
$L/\Lambda_1$ is called an {\em affine automorphism\/} of
$L/\Lambda_1$ if there exists $\sigma\in\Aut(L)$ such that
$\tau(gx)=\sigma(g)\tau(x)$ for all $x\in L/\Lambda_1$. The group of
all affine automorphisms of $L/\Lambda_1$ is denoted by
$\Aff(L/\Lambda_1)$. It is endowed with the compact-open topology;
i.e.\ its open subbase consists of the sets of the form
$\{\tau\in\Aff(L/\Lambda_1):\tau(C)\subset U\}$, where $C$ is a
compact subset of $L/\Lambda_1$ and $U$ is an open subset of
$L/\Lambda_1$. 
\end{definition}

\begin{remark}  \label{rem:affine} \rm
(1) $\Aff(L/\Lambda_1)$ is a locally compact topological group acting
continuously on $L/\Lambda_1$.  (2) If $\sigma\in\Aut(L)$ is such that
$\sigma(\Lambda_1)=\Lambda_1$ and if $g\in L$, then the map $\tau$ on
$L/\Lambda_1$ defined by $\tau(h\Lambda_1)=g\sigma(h)\Lambda_1$ for
all $h\in L$ is an affine automorphism. (3) Let
$\Lambda_1'$ be the maximal closed normal subgroup of
$L$ contained in $\Lambda_1$. Define $\bar L=L/\Lambda'_1$ and
$\bar\Lambda_1=\Lambda_1/\Lambda'_1$. Then we have natural
isomorphisms $L/\Lambda\cong \bar L/\bar\Lambda_1$ and
$\Aff(L/\Lambda_1)=\Aff(\bar L/\bar\Lambda_1)$.
\end{remark}

\begin{theorem} \label{thm:factors}
Let $G$ be a subgroup of $L$ which is generated by one-parameter
unipotent subgroups of $L$ contained in $G$. Suppose that $G$ acts
ergodically on $L/\Lambda$. Let $X$ be a Hausdorff locally compact
space with a continuous $G$-action and $\phi:L/\Lambda\to X$ a
continuous surjective $G$-equivariant map. Then there exists a closed
subgroup $\Lambda_1$ containing $\Lambda$, a compact group $ K$
contained in the centralizer of the subgroup of translations by
elements of $G$ in $\Aff(L/\Lambda_1)$, and a $G$-equivariant
continuous surjective map $\psi: K\bs L/\Lambda_1\to X$ such that the
following holds:

\begin{enumerate}
\item
If $\rho:L/\Lambda\to K\bs L/\Lambda_1$ is defined by $\rho(g\Lambda)=
K(g\Lambda_1)$, $\forall\, g\in L$, then $\rho$ is $G$-equivariant and
\[\phi=\psi\circ\rho.\]
\item 
Given a neighbourhood $\Omega$ of $e$ in $Z_L(G)$, there exists an open
dense $G$-invariant subset $X_0$ of $L/\Lambda_1$ such that for any
$x\in X_0$ and $y\in L/\Lambda_1$ if $\psi(K(x))=\psi(K(y))$ then
$y\in K(\Omega x)$. In this situation, further if $\cl{Gx}=L/\Lambda_1$,
then $ K(y)= K(x)$.  
\end{enumerate}
\end{theorem}

The above description of topological factors of unipotent flows is
also of independent interest. The measurable factors of unipotent
flows were described by Witte~\cite{Witte}. 

The next result is an immediate consequence of
theorems~\ref{thm:G/P}~and~\ref{thm:factors}.

\begin{corollary} \label{cor:topological:factors} Let $L$ be a Lie group,
$\Lambda$ a lattice in $L$, and $G$ a connected semisimple Lie group
with finite center, realized as a closed subgroup of $L$.  Suppose
that the action of $G_1$ on $L/\Lambda$ is minimal for any closed
normal subgroup $G_1$ of $G$ such that $\RR$-$\rank(G/G_1)\leq 1$. Let
$Y$ be a locally compact Hausdorff space with a continuous $G$-action,
$P$ a parabolic subgroup of $G$, and $\phi:L/\Lambda\times G/P\to Y$ a
continuous surjective $G$-equivariant map, where $G$-acts diagonally
on $L/\Lambda\times G/P$. Then there exist a parabolic subgroup $Q$ of
$G$ containing $P$, a closed subgroup $\Lambda_1$ of $L$ containing
$\Lambda$, and a compact group $ K$ contained in the centralizer of
the image of $G$ in $\Aff(L/\Lambda_1)$, such that $Y$ is
$G$-equivariantly homeomorphic to $( K\backslash L/\Lambda_1)\times
(G/Q)$ and $\phi$ is the natural quotient map.

In particular if, as in question~\ref{que:stuck}, there exists a map
$\psi:Y\to L/\Lambda$ such that $\psi\circ\phi$ is the projection on the
first factor, then $\Lambda_1=\Lambda$ and $ K$ is trivial. Hence
$Y$ is $G$-equivariantly homeomorphic to $L/\Lambda\times G/Q$ and
$\phi$ is the natural quotient map. 
\end{corollary}

For the purpose of other applications, we obtain stronger a measure
theoretic version of theorem~\ref{thm:parabolic}. Before the
statement, we recall some definitions.

For any Borel map $T:X\to Y$ of Borel spaces and a Borel measure
$\lambda$ on $X$, the Borel measure $T_\ast \lambda$ defined by
$T_\ast\lambda(E)=\lambda(T\inv(E))$, for all Borel sets $E\subset Y$,
is called the image of $\lambda$ under $T$.

For any Borel measure $\mu$ on $L/\Lambda$ and any $g\in L$, the
translated measure $g\cdot\mu$ on $L/\Lambda$ is the image of $\mu$ under
the map $x\mapsto gx$ on $L/\Lambda$. 

On a locally compact space $X$, for a sequence $\{\mu_i\}$ of finite
Borel measures and $\mu$ a finite Borel measure, we say that
$\mu_i\to\mu$ as $i\to\infty$, if and only if for all bounded
continuous function $f$ on $X$, $\scriptsize{\int_X f\,d\mu_i\to\int_X
f\,d\mu}$ as $i\to\infty$.

\begin{notation}   \label{not:Delta} \rm
Let $G$ be a connected semisimple real algebraic group. Let $A$ be an
$\RR$-split torus in $G$ such that the set of real roots on $A$ for
the adjoint action on the Lie algebra of $G$ forms a root system. Fix
an order on this set of roots and let $\Delta$ be the corresponding
system of simple roots. Let $\cl{A}^+$ be the closure of the positive
Weyl chamber in $A$. Let $\{a_i\}_{i\in\NN}$ be a sequence in
$\cl{A}^+$ such that for any $\alpha\in\Delta$, either
$\sup_{i\in\NN}\alpha(a_i)<\infty$ or $\alpha(a_i)\to\infty$ as
$i\to\infty$. Put
\[
U^+=\{g\in G:a_i\inv g a_i \to e\mbox{ as $i\to\infty$}\}.
\]
\end{notation}

\begin{theorem} \label{thm:translate:horosphere}
Consider the notation~\ref{not:Delta}.  Assume that $U^+$ is not
contained in any proper closed normal subgroup of $G$.  Suppose that
$G$ is realized as a Lie subgroup of $L$ and that $\pi(G)$ is dense in
$L/\Lambda$. Then for any probability measure $\lambda$ on $U^+$ which
is absolutely continuous with respect to a Haar measure on $U^+$,
\[
a_i\cdot\pi_\ast(\lambda)\to\mu_L, \qquad \mbox{as $i\to\infty$}.
\]
In other words, for any bounded continuous function
$f$ on $L/\Lambda$, 
\[
\lim_{i\to\infty}\int_{U^+} f(a_i\pi(\omega))\,d\lambda(\omega) =
\int_{L/\Lambda}f\,d\mu_L.
\]

In particular, for any Borel set $\Omega$ of $U^+$ having
strictly positive Haar measure,
\[
\cl{\bigcup_{i\in\NN}a_i\cdot\pi(\Omega)}=L/\Lambda.
\]
\end{theorem}

Using this theorem we obtain the following generalization of a result
due to Duke, Rudnick and Sarnak~\cite{DRS:affine}; their result
corresponds to the case of $L=G$. First we need a definition.

Let $G$ be a semisimple Lie group. A subgroup $S$ of $G$ is said to be
{\em symmetric\/} if there exists an involution $\sigma$ of $G$
(i.e. $\sigma$ is a continuous automorphism and $\sigma^2=1$) such
that $S=\{g\in G:\sigma(g)=g\}$. For example, any maximal compact
subgroup of $G$ is a symmetric subgroup, for it is the fixed point set
of a Cartan involution of $G$.

\begin{corollary} \label{cor:symmetric}
Let $G$ be a connected real algebraic semisimple Lie group realized as
a Lie subgroup of $L$, $S$ the connected component of the identity of
a symmetric subgroup of $G$, and $\{g_i\}_{i\in\NN}$ a sequence
contained in $G$.  Suppose that $\pi(S)$ is closed and admits an
$S$-invariant probability measure, say $\mu_S$. Also suppose that
$\pi(G_1)$ is dense in $L/\Lambda$, for any closed normal subgroup
$G_1$ of $G$ such that the image of $\{g_i\}$ in $G/(SG_1)$ admits a
convergent subsequence. Then the sequence of measures $g_i\cdot\mu_S$
converges to $\mu_L$; that is, for every bounded continuous function
$f$ on $L/\Lambda$,
\[
\lim_{i\to\infty} \int_{\pi(S)} f(g_ix)\,d\mu_S(x) = \int_{L/\Lambda}
f\,d\mu_L. 
\]
\end{corollary}

In the case of $L=G$, Eskin and McMullen~\cite{EM:mixing} gave a proof of
this result using the mixing property of geodesic flows. The main
technical observation in their proof is what they call `a wave front
lemma'. In the general case of $L\supset G$, our analogue of the wave
front lemma is theorem~\ref{thm:translate:horosphere}. 

Using the arguments of the proof of corollary~\ref{cor:symmetric}, one
can also deduce the following result from
theorem~\ref{thm:translate:horosphere}. 

\begin{corollary}  \label{cor:mixing?}
Let $G$ be a connected real algebraic semisimple group realized as a
Lie subgroup of $L$. Let $\{g_i\}$ be a sequence in $G$. Suppose that
$\pi(G_1)$ is dense in $L/\Lambda$ for any closed normal subgroup
$G_1$ of $G$ such that the image of $\{g_i\}$ in $G/G_1$ admits a
convergent subsequence. Then for any Borel probability measure
$\lambda$ on $G$ which is absolutely continuous with respect to a Haar
measure on $G$, 
\[
g_i\pi_\ast(\lambda)\to\mu_L \qquad \mbox{as $i\to\infty$}.
\]
In particular, for any Borel set $\Omega$ of $G$ having strictly
positive Haar measure,
\[
\cl{\bigcup_{i\in\NN}g_i\cdot\pi(\Omega)}=L/\Lambda.
\]
\end{corollary} 

Using the method of our proof, one can also obtain the uniform
versions of theorem~\ref{thm:translate:horosphere} and
corollary~\ref{cor:symmetric} which are similar to 
\cite[Theorems~4.3-4]{EMM:upperbound}.

The the main result of this paper is
theorem~\ref{thm:translate:horosphere} and other results (except
theorem~\ref{thm:factors}) are derived from it. The main steps of its
proof are as follows.  First suppose that the set of probability
measures $\{a_i\cdot\pi_\ast(\lambda):i\in\NN\}$ is not is not
relatively compact in the space of all probability measures on
$L/\Lambda$. Using an extension of a result of Dani and
Margulis~\cite{DM:asymptotic}, in section~2 we see that there exist a
nonempty open set $\Omega\subset U^+$, a finite dimensional
representation $V$ of $L$, a discrete set $\{\bfv_i:i\in\NN\}\subset
V$, and a compact set $K\subset V$ such that that
$a_i\Omega\cdot\bfv_i\subset K$ for infinitely many $i\in\NN$. Via
some elementary observations about representations of semisimple Lie
groups, in section~5 we show that the conditions mentioned above lead
to a contradiction when we restrict the representation to $G$. Now let
a probability measure $\mu$ be a limit distribution of the sequence
$\{a_i\cdot\pi_\ast(\lambda)\}$. We observe that $\mu$ is
$U^+$-invariant. Using Ratner's~\cite{R:measure} description of finite
measures on $L/\Lambda$ which are ergodic and invariant under the
action of a unipotent subgroup, in section~3 we conclude that either
$\mu=\mu_L$, or $\mu$ is nonzero when restricted to the image under
$\pi$ of some strictly lower dimensional `algebraic subvariety' of
$L$. Using techniques developed in
\cite{DM:generic,Sh:equidistribution,DM:equidistribution,MS:limit}, in
section~4 we see that in the later case the above type of condition on
a finite dimensional representation of $L$ must hold, and this again
leads to a contradiction. Thus $\mu=\mu_L$ and hence $\mu_L$ is the
only limit distribution of $\{a_i\cdot\mu_\Omega\}$.

\section{A condition for returning to compact sets}

In~\cite{DM:asymptotic} Dani and Margulis proved that large compact
sets in finite volume homogeneous spaces have relative measures close
to $1$ on the trajectories of unipotent flows which originate from a
fixed compact set. This result was generalized in \cite{EMS:nonzero}
to a larger class of higher dimensional trajectories. In these results
one considered only the case of arithmetic lattices in algebraic
semisimple Lie groups defined over $\QQ$. Here we modify them to
include the case of any lattice in any Lie group.

\begin{notation} \rm
Let $G$ be a Lie group and $\la{G}$ the Lie algebra associated to $G$. For
$d,m\in\NN$, let $\cP_{d,m}(G)$ denote the set of continuous maps
$\Theta:\RR^m\to G$ such that for all $\bfc,\bfa\in\RR^m$ and
$X\in\la{G}$, the map
\[
t\in\RR\mapsto \Ad\circ\Theta(t\bfc+\bfa)(X)\in\la{G}
\]
is a polynomial of degree at most $d$ in each co-ordinate of $\la{G}$
(with respect to any basis).

We shall write  $\cP_d(G)$ for the set $\cP_{d,1}(G)$. 
\end{notation}

\begin{theorem}[Dani, Margulis] \label{thm:compact:compact}
Let $G$ be a Lie group, $\Gamma$ a lattice in $G$, and $\pi:G\to G/\Gamma$
the natural quotient map. Then given a
compact set $C\subset G/\Gamma$, an $\epsilon>0$, and a $d\in\NN$, there
exists a compact subset $K \subset G/\Gamma$ with the following property:
For any $\Theta\in\cP_{d,m}(G)$ and any bounded open
convex set $B\subset\RR^m$, one of the following conditions hold:
\begin{enumerate}
\item
$(1/\nu(B))\nu(\{\bft\in B:\pi(\Theta(\bft)) \in K\})\geq
1-\epsilon$, where $\nu$ denotes the Lebesgue measure on $\RR^m$.
\item
$\pi(\Theta(B))\cap C = \emptyset$.
\end{enumerate}
\end{theorem}

\proof
See~\cite[Theorem~3.1]{Shah:polynomial}.
\qed

The usefulness of the above result is enhanced by the following
theorem which provides an algebraic condition in place of the
geometric condition $\pi(\Theta(B))\cap C=\emptyset$.

\begin{notation} \label{not:V} \rm
Let $G$ be a connected Lie group and $\la{G}$ denote the
Lie algebra associated to $G$. Let
$V_G=\oplus_{k=1}^{\dim\la{G}}\wedge^k\la{G}$, the direct sum of exterior
powers of $\la{G}$, and consider the linear $G$-action on $V_G$ via the
representation $\oplus_{l=1}^{\dim\la{G}}\wedge^l\Ad$, the direct sum of
exterior powers of the adjoint representation of $G$ on $\la{G}$.

Fix any Euclidean norm on $\la{G}$ and let
$\cB=\{\bfe_1,\ldots,\bfe_{\dim\la{G}}\}$ denote an orthonormal basis of
$\la{G}$. There is a unique Euclidean norm $\|\cdot\|$ on $V_G$ such
that the associated basis of $V_G$ given by
\[
\{\bfe_{l_1}\wedge\cdots\wedge \bfe_{l_r}:1 \leq l_1 <\ldots< \l_r \leq
\dim\la{G},\, r=1,\ldots,\dim\la{G}\}
\]
is orthonormal. This norm is independent of the choice of $\cB$. 

To any Lie subgroup $W$ of $G$ and the associated Lie subalgebra $\la{W}$
of $\la{G}$ we associate a unit-norm vector
$\bfp_W\in\wedge^{\dim\la{W}}\la{W}\in V_G$. 
\end{notation}

\begin{theorem}[Cf.~\cite{DM:asymptotic}]  \label{thm:return:strong}
Let $G$ be a connected Lie group, $\Gamma$ a lattice in $G$, and
$\pi:G\to\GmG$ the natural quotient map.  Let $M$ be the smallest
closed normal subgroup of $G$ such that $\bar G=G/M$ is a semisimple
group with trivial center and without nontrivial compact normal
subgroups. Let $q:G\to\bar G$ be the quotient homomorphism. Then there
exist finitely many closed subgroups $W_1,\ldots,W_r$ of $G$ such that
each $W_i$ is of the form $q\inv(U_i)$ with $U_i$ the unipotent
radical of a maximal parabolic subgroup of $\bar G$, $\pi(W_i)$ is
compact and the following holds: Given $d,m\in\NN$ and reals
$\alpha,\epsilon>0$, there exists a compact set $C\subset \GmG$ such
that for any $\Theta\in\cP_{d,m}(G)$, and a bounded open convex set
$B\subset\RR^m$, one of the following conditions is satisfied:
\begin{enumerate}
\item 
There exist $\gamma\in\Gamma$ and $i\in\{1,\ldots,r\}$ such that 
\[
\sup_{\bft\in B}\|\Theta(\bft)\gamma\cdot\bfp_{W_i}\|<\alpha.
\]

\item
$\pi(\Theta(B))\cap C\neq\emptyset$, and hence condition~(1) of
theorem~\ref{thm:compact:compact} holds.
\end{enumerate}
\end{theorem}

\proof Let $R$ be the radical of $G$, $C$ the maximal connected
compact normal subgroup of $G/R$, $S=(G/R)/C$ and $Z$ the center of
$S$. Note that $S$ is a semisimple Lie group without proper compact
connected normal subgroups. Clearly $S/Z\cong G/M$. Therefore $M$ is
the inverse image of $Z$ in $G$.

Let $H=\cl{R\Gamma}^0$. Then $H\Gamma$ is closed and $H\cap\Gamma$ is
a lattice in $H$ (see~\cite[Lemma~1.7]{Rag:book}). By Auslander's
theorem~\cite[Theorem~8.24]{Rag:book} $H$ is solvable, and so is
its image in $S$. By Borel's density theorem~\cite[Lemma~5.4,
Corollary~5.16]{Rag:book} the image is a normal subgroup of $S$ and
therefore it has to be trivial. Hence $H\subset M^0$, and since
$R\subset H$, $M^0/H$ is compact. Since $H$ is solvable, by Mostow's
theorem~\cite[Theorem~3.1]{Rag:book} $H/(H\cap\Gamma)$ is
compact. Therefore $M^0/H\cap\Gamma$ is compact. So $M^0\Gamma/\Gamma$
is compact and $M^0\Gamma$ is closed.  

Therefore the image $\Delta$ of $\Gamma$ in $S$ is discrete, and hence
a lattice in $S$. Therefore by Borel's density
theorem~\cite[Corollary~5.18]{Rag:book} $Z\Delta$ is discrete. Hence
$\Delta$ is of finite index in $Z\Delta$ and hence $M^0\Gamma$ is of
finite index in $M\Gamma$. Hence $M\Gamma/\Gamma$ is compact, i.e.\
$\pi(M)$ is compact.

Thus $\bar\Gamma=q(\Gamma)$ is a lattice in $\bar G$ and the fibers of
the map $\bar q:G/\Gamma\to\bar G/\bar\Gamma$ are compact
$M$-orbits. Therefore without loss of
generality, we may assume that $\bar G=G$.

\clearpage

Then there are finitely many normal connected subgroups
$G_1,\ldots,G_r$ of $G$ such that $G=G_1\times\cdots\times G_r$ and
each $\Gamma_i=G_i\cap\Gamma$ is an irreducible lattice in $G_i$
(see~\cite[Sect.~5.22]{Rag:book}). In proving the theorem without loss
of generality we may replace $\Gamma$ by its finite-index subgroup
$\Gamma_1\times\cdots\times\Gamma_r$. In order to prove the theorem
for $G$, it is enough to prove it for each $G_i$ separately.  Thus
without loss of generality we may assume that $\Gamma$ is an
irreducible lattice.

The result in the case of $\RR$-rank$(G)=1$ can be deduced from the
arguments  in~\cite[(2.4)]{Dani:rk1}. 

Next suppose that $\RR\mbox{-rank}(G)\geq 2$. Then by the
arithmeticity theorem of Margulis~\cite{Mar:arith}, $\Gamma$ is
an arithmetic lattice. Therefore there exist a semisimple algebraic
group $\bfG$ defined over $\QQ$ and a surjective homomorphism
$\rho:\bfG(\RR)^0\to G$ with compact kernel such that, for
$\Lambda=\bfG(\ZZ)\cap\bfG(\RR)^0$, the subgroup
$\Gamma\cap\rho(\Lambda)$ is a subgroup of finite index in both
$\Gamma$ and $\rho(\Lambda)$. Again without loss of
generality we may replace $G$ by $\bfG(\RR)^0$ and $\Gamma$ by
$\Lambda$.  In this case the result follows from
\cite[Thm.~3.6]{EMS:nonzero}.  \qed


\section{Description of measures invariant under a unipotent flow}

In this and the next section, let $G$ denote a Lie group, $\Gamma$ a
lattice in $G$, and $\pi:G\to G/\Gamma$ the natural quotient map.

A subgroup $U$ of $G$ is called {\em unipotent\/} if $\Ad u$ is a
unipotent endomorphism of the Lie algebra of $G$ for every $u\in U$. 

Let $\cH_\Gamma$ denote the collection of all closed connected
subgroups $H$ of $G$ such that (1) $H\supset\Gamma$, (2)
$H/H\cap\Gamma$ admits a finite $H$-invariant measure, and (3) the
subgroup generated by all one-parameter unipotent subgroups of $H$
acts ergodically on $H/H\cap\Gamma$ with respect to the $H$-invariant
probability measure.  In particular, the Zariski closure of
$\Ad(H\cap\Gamma)$ contains $\Ad(H)$
(see~\cite[Theorem~2.3]{Sh:equidistribution}).

\newcommand{\citfive}{{\rm \cite[Theorem~1.1]{R:measure}}}
\begin{theorem}[\citfive] \label{thm:countable}
The collection $\cH_\Gamma$ is countable.
\end{theorem}

Let $W$ be a subgroup of $G$ which is generated by one-parameter
unipotent subgroups of $G$ contained in $W$. For any $H\in\cH_\Gamma$,
define
\begin{eqnarray*}
N_G(H,W)&=&\{g\in G: W\subset gHg\inv\}, \\[3pt]
S_G(H,W)&=&\bigcup_{\stackrel{H'\in\cH_\Gamma,\,H'\subset H}{\dim
H'<\dim H}} N_G(H',W).
\end{eqnarray*}

Note that (see \cite[Lemma~2.4]{MS:limit}),
\begin{equation}  \label{eq:pi(N(H,W)-S(H,W))}
\pi(N_G(H,W)\setminus S_G(H,W))=\pi(N_G(H,W)) \setminus \pi(S_G(H,W)).
\end{equation}

We reformulate Ratner's classification~\cite{R:measure} of finite
measures which are invariant and ergodic under unipotent flows on
homogeneous spaces of Lie groups, using the above definitions
(see~\cite[Theorem 2.2]{MS:limit}).

\begin{theorem} \label{thm:mu_H}
Let $W$ be a subgroup as above and $\mu$ a $W$-invariant probability
measure on $G/\Gamma$. For every $H\in\cH_\Gamma$, let $\mu_H$ denote
the restriction of $\mu$ on $\pi(N_G(H,W)\setminus S_G(H,W))$. Then
the following holds.
\begin{enumerate}
\item
The measure $\mu_H$ is $W$-invariant, and any $W$-ergodic component
of $\mu_H$ is of the form $g\cdot\lambda$, where  $g\in N_G(H,W)\setminus
S_G(H,W)$ and $\lambda$ is a $H$-invariant measure on $H\Gamma/\Gamma$. 
\item
For any Borel measurable set $A\subset G/\Gamma$,
\[
\mu(A)=\sum_{H\in\cH_\Gamma^\ast}\mu_H(A),
\]
where $\cH_\Gamma^\ast\subset\cH_\Gamma$ is a countable set consisting of one
representative from each $\Gamma$-conjugacy class of elements in $\cH_\Gamma$.
\end{enumerate}

In particular, if $\mu(\pi(S(G,W))=0$ then $\mu$ is the unique
$G$-invariant probability measure on $G/\Gamma$.
\end{theorem}


\section{Linear presentation of $G$-actions near singular sets}

\label{sec:linearise}

Let $C\subset\pi(N_G(H,W)\setminus S_G(H,W))$ be any compact set.  It
turns out that on certain neighborhoods of $C$ in $G/\Gamma$, the
$G$-action is equivariant with the linear $G$-action on certain
neighbourhoods of a compact subset of a linear subspace in a finite
dimensional linear $G$-space. We study unipotent trajectories in those
thin neighbourhoods of $C$ via this linearisation. This type of
technique is developed in
(\cite{DM:generic,Sh:equidistribution,DM:equidistribution,Shah:polynomial,MS:limit,EMS:counting}).

Let $V_G$ be the representation of $G$ as described in
notation~\ref{not:V}.  For $H\in\cH_\Gamma$, let $\eta_H:G\to V_G$ be
the map defined by $\eta_H(g)=g\bfp_H=(\wedge^d\Ad g)\bfp_H$ for all
$g\in G$. Let $N_G(H)$ denotes the normaliser of $H$ in $G$. Define
\[
N_G^1(H)=\eta_H\inv(\bfp_H)=\{g\in N_G(H) : \det(\Ad g|_\la{H})=1\}.
\]

\begin{proposition} {\rm(\cite[Theorem~3.4]{DM:equidistribution})}  
\label{prop:discrete}
The orbit $\Gamma\cdot\bfp_H$ is closed, and hence discrete. In
particular, the orbit $N_G^1(H)\Gamma/\Gamma$ is closed in $G/\Gamma$. 
\end{proposition}

Let $W$ be a subgroup which is generated by one-parameter unipotent
subgroups of $G$ contained in $W$.  

\begin{proposition}  {\rm (\cite[Prop.~3.2]{DM:equidistribution})} 
\label{prop:variety}
Let $V_G(H,W)$ denote the linear span of $\eta(N_G(H,W))$ in $V_G$. Then
\[
\eta_H\inv(V_G(H,W))=N_G(H,W).
\]
\end{proposition}

\begin{theorem} \label{thm:general-avoid}
Let $\epsilon>0$, $d,m\in\NN$, and a compact set $C\subset
\pi(N_G(H,W)\setminus S_G(H,W))$ be given. Then there exists a compact set
$D\subset V_G(H,W)$ such that given any neighbourhood $\Phi$ of $D$ in
$V_G$, there exists a neighbourhood $\Psi$ of $C$ in $G/\Gamma$ such that
for any $\Theta\in\cP_{d,m}(G)$, and a bounded open convex set $B\subset
\RR^m$, one of the following conditions is satisfied:
\begin{enumerate}
\item 
$\Theta(B)\gamma\cdot\bfp_H \subset \Phi$ for some $\gamma\in\Gamma$. 
\item 
\[
\frac{1}{\nu(B)}\nu(\{\bft\in B:\Theta(\bft)\Gamma/\Gamma\in\Psi\})
<\epsilon. 
\]
\end{enumerate}
\end{theorem}

\proof 
The result is easily deduced from \cite[Prop.~5.4]{Shah:polynomial}.
See also the proof of \cite[Thm.~5.2]{Shah:polynomial}.
\qed

\subsection*{Some related results on unipotent flows}

We recall a theorem of Ratner~\cite{R:equidistribution} on closures
of individual orbits of unipotent flows.

\begin{theorem}[Ratner]  \label{thm:Ratner}
Let $G$, $\Gamma$ and $W$ be as above. Then for any $x\in G/\Gamma$, there
exists a closed subgroup $F$ of $G$ containing $W$ such that $\cl{Wx}=Fx$
and the orbit $Fx$ admits a unique $F$-invariant probability measure, say
$\mu_F$. Also $\mu_F$ is $W$-ergodic. 
\end{theorem}

Next we recall a version of the main result of~\cite{MS:limit}.

\begin{theorem}[\cite{MS:limit}] \label{thm:limit}
Let $x\in G/\Gamma$, and sequences $\{F_i\}$ of closed subgroups of
$G$ and $g_i\to e$ in $G$ be such that each of the orbits
$F_i(g_ix)$ is closed, and admits an $F_i$-invariant probability
measure, say $\mu_i$. Suppose that the subgroup generated by all
unipotent one-parameter subgroups of $G$ contained in $F_i$ acts
ergodically with respect to $\mu_i$, $\forall i\in\NN$. Then there
exits a closed subgroup $F$ of $G$ such that the orbit $Fx$ is closed,
and admits a $F$-invariant probability measure, say $\mu$, and a
subsequence of $\{\mu_i\}$ converges to $\mu$.

Moreover if $\mu_i\to\mu$ as $i\to\infty$, then $g_i\inv F_i
g_i\subset F$ for all large $i\in\NN$.
\end{theorem}

\section{Some results on linear representations}

In view of the proposition~\ref{prop:discrete}, in order to obtain
further consequences when either condition~1 of
theorem~\ref{thm:return:strong} or condition~1 of
theorem~\ref{thm:general-avoid} holds for a sequence
$\{\Theta_i\}\subset\cP_{d,m}(G)$, the following elementary result is
very useful.

\subsection*{Linear actions of Unipotent subgroups}

\begin{lemma}  \label{lemma:nilpotent}
Let $V$ be a finite dimensional real vector space equipped with a
Euclidean norm. Let\/ $\la{N}$ be a nilpotent Lie subalgebra of
$\End(V)$. Let $N$ be the associated unipotent subgroup of
$\Aut(V)$. Let $\{\bfb_1,\ldots,\bfb_m\}$ be a basis of
$\la{N}$. Consider the map $\Theta:\RR^m\to N$ defined as
\[
\Theta(t_1,\ldots,t_m)=\exp(t_m\bfb_m)\cdots\exp(t_1\bfb_1),
\qquad \forall(t_1,\ldots,t_m)\in\RR^m.
\]
For $\rho>0$, define
\[
B_\rho=\{\Theta(t_1,\ldots,t_m)\in N:0\leq t_k<\rho,\, k=1,\ldots,m\}.
\]
Put 
\[
W=\{\bfv\in V: n\cdot\bfv=\bfv,\,\forall n\in N\}.
\]
Let $\pr_W$ denote the orthogonal projection on $W$. Then for any
$\rho>0$, there exists $c>0$ such that for every $\bfv\in V$,
\[
\|\bfv\|\leq c\cdot\sup_{\bft\in B_\rho}\|\pr_W(\Theta(\bft)\cdot\bfv)\|.
\]
\end{lemma}

\proof
For $k=1,\ldots,m$, let $n_k\in\NN$ be such that $\bfb_k^{n_k}=0$. Let 
\[
\cI=\{I=(i_1,\ldots,i_m):0\leq i_k\leq n_k-1,\,k=1,\ldots,m\}.
\]
For $\bft=(t_1,\ldots,t_m)\in\RR^m$ and $I=(i_1,\ldots,i_m)\in\cI$, define
\[
\bft^I=t_m^{i_m}\cdots t_1^{i_1} \qquad \mbox{and}\qquad 
\bfb^I=\frac{\bfb_m^{i_m}\cdots\bfb_1^{i_1}}{i_m!\cdots i_1!}.
\]

Then for all $\bfv\in V$ and $\bft\in\RR^m$, we have
\begin{equation} \label{eq:Theta(t)}
\Theta(\bft)\cdot \bfv=\sum_{I\in\cI}\bft^I\cdot(\bfb^I \bfv).
\end{equation}

We define a transformation $T:V\to\oplus_{I\in\cI}W$ by
\begin{equation} \label{eq:T(v)}
T(\bfv)=\left(\pr_W(\bfb^I\cdot\bfv)\right)_{I\in\cI}, \qquad \forall \bfv\in V.
\end{equation}
We claim that $T$ is injective. To see this, suppose there exists $\bfv\in
V\setminus\{0\}$ such that $T(\bfv)=0$. Then $N\cdot\bfv\subset W^\perp$,
the orthogonal complement of $W$. Hence $W^\perp$ contains a nontrivial
$N$-invariant subspace. Then by Lie-Kolchin theorem, $W^\perp$ contains a
nonzero vector fixed by $N$. Then $W\cap W^\perp\neq\{0\}$, which is a
contradiction. 

We consider $\oplus_{I\in\cI}V$ equipped with a box norm; that is
\[
\|(v_I)_{I\in\cI}\|=\sup_{I\in\cI}\|v_I\|,\qquad  
\mbox{where $v_I\in V$, $\forall I\in\cI$}.
\]
Since $T$ is injective, there exists a constant $c_1>0$ such that
\[
\|\bfv\|\leq c_1\cdot\|T(\bfv)\|, \qquad \forall \bfv\in V.
\]

For all $k=1,\ldots,m$, and $j_k=1,\ldots,n_k$, fix
$0<t_{k,1}<\cdots<t_{k,n_k}<\rho$ and put 
$M_k=\left(t_{k,j_k}^{i_k}\right)_{0\leq i_k\leq n_k-1,\, 1\leq j_k\leq
n_k}$ for $k=1,\ldots,m$. Then $\det M_k$ is a Vandermonde determinant and
hence $M_k$ is invertible. 

Let 
\[
\cJ=\{J=(j_1,\ldots,j_m):1\leq j_k\leq n_k,\,k=1,\ldots,m\}.
\]
For $J=(j_1,\ldots,j_m)\in\cJ$, put 
\[
\bft_J=(t_{1,j_1},\ldots,t_{m,j_m})\qquad 
\mbox{and}\qquad M=\left(\bft_J^I\right)_{(I,J)\in\cI\times\cJ}.
\]

Take $\bfv\in V$. Put
\[ 
X_\cI=T(\bfv)\qquad \mbox{and}\qquad 
Y_\cJ=\left(\pr_W(\Theta(\bft_J)\bfv)\right)_{J\in\cJ}. 
\]
Then by equations~\ref{eq:Theta(t)} and \ref{eq:T(v)},
\[
M\cdot X_\cI=Y_\cJ.
\]
Since $M=M_1\otimes\cdots\otimes M_m$ and each $M_k$ is invertible, we
have that $M$ is invertible. Hence  
\[
X_\cI=M\inv\cdot Y_\cJ.
\]
Put $c_2=\|M\inv\|$ and $c=c_1c_2$. Then
\[
\|\bfv\|\leq c_1\|T(\bfv)\|=c_1\|X_\cI\|\leq c_1c_2\|Y_\cJ\|=
c\cdot\sup_{J\in\cJ}\|\pr_W(\Theta(\bft_J)\bfv)\|. 
\]
This completes the proof. 
\qed

\subsection*{Linear actions of semisimple groups}

We fix the following setup for the rest of this section.

\begin{notation} \label{not:normal} \rm 
Consider the notation~\ref{not:Delta}.  Put
\[
\Phi=\{\alpha\in\Delta:\alpha(a_i)\to\infty\mbox{ as } i\to\infty\}.
\]
Let $P^+$ be the standard parabolic subgroup associated to the set of
roots $\Delta\setminus\Phi$. Then $U^+=\{g\in G:a_i\inv g a_i\to e$ as
$i\to\infty\}$ is the unipotent radical of $P^+$. Let $P^-$ denote the
standard opposite parabolic subgroup for $P^+$ and let $U^-$ be the
unipotent radical of $P^-$. Note that
\begin{equation} \label{eq:relcpt}
P^-=\{g\in G: \cl{\{a_iga_i\inv:i\in\NN\}}\mbox{ is compact}\}.
\end{equation}
Also put $Z=P^-\cap P^+$. Then $P^-=U^-Z$. Let $\la{G}$, $\la{U}^-$,
$\la{Z}$, and $\la{U}^+$ denote the Lie algebras associated to $G$,
$U^-$, $Z$, and $U^+$, respectively. Then
\begin{equation}  \label{eq:G=U-ZU+}
\la{G}=\la{U}^-\oplus \la{Z}\oplus \la{U}^+.
\end{equation}
\end{notation}

\begin{lemma} \label{lemma:take:far}
Consider a continuous nontrivial irreducible representation of $G$ on
a finite dimensional normed vector space $V$. Let $W=\{\bfv\in
V:W\cdot\bfv=\bfv\}$. Let $\{\bfv_i\}\subset W$ be a sequence such
that $\inf_{i\in\NN}\|\bfv_i\|>0$. Then
\[
\|a_i\cdot\bfv_i\|\to\infty\qquad 
\mbox{as $i\to\infty$}.
\]
\end{lemma}

\proof 
Since $A$ is $\RR$-split, there is a finite set $\Lambda$ of real
characters on $A$ such that for each $\lambda\in\Lambda$, if we define
\[
V_\lambda=\{\bfv\in V:a\cdot\bfv=\lambda(a)\bfv,\,\forall a\in A\},
\]
then $V=\oplus_{\lambda\in\Lambda} V_\lambda$. After passing to an appropriate
subsequence, if we define 
\begin{eqnarray*}
\Lambda_+&=&\{\lambda\in\Lambda:\lambda(a_i)\to\infty 
\mbox{ as $i\to\infty$}\}\\ 
\Lambda_-&=&\{\lambda\in\Lambda:\lambda(a_i)\to 0 
\mbox{ as $i\to\infty$}\}, \qquad \mbox{and }\\
\Lambda_0&=&\{\lambda\in\Lambda:\lambda(a_i)\to c 
\mbox{ for some $c>0$ as $i\to\infty$}\},
\end{eqnarray*}
then $\Lambda=\Lambda_+\cup \Lambda_0 \cup \Lambda_-$.

Since $U^+$ is normalized by $A$, we have that $W$ is invariant under
the action of $A$.  Therefore
\[
W=\oplus_{\lambda\in\Lambda} (W\cap V_\lambda).
\]

Suppose that there exists $\bfw\in W\cap V_\lambda\setminus\{0\}$ for
some $\lambda\in\Lambda_0\cup\Lambda_-$.  For any $g\in P^-$, we have
$a_iga_i\inv\to g_0$ for some $g_0\in P^-$.  Therefore as
$i\to\infty$,
\[
a_i(g\bfw)=a_iga_i\inv(a_i\bfw)\to c(g_0\bfw) \mbox{\qquad for some
$c\geq 0$.}
\]
Hence
$P^-\bfw\subset\oplus_{\lambda\in\Lambda_0\cup\Lambda_-}V_\lambda$.
Now $U^+\bfw=\bfw$ and by notation~\ref{not:normal} $P^-U^+$ is open
in $G$.  Therefore
$G\cdot\bfw\subset\oplus_{\lambda\in\Lambda_0\cup\Lambda_-}V_\lambda$.
Since $V$ is irreducible, $\Lambda=\Lambda_0\cup\Lambda_-$. Now since
$G$ is semisimple, $\det g=1$ for all $g\in G$ and hence
$\Lambda_-=\emptyset$.  Thus $\Lambda=\Lambda_0$. 

Now for any relatively compact neighbourhood $\Omega$ of $U^+$ and any
$\bfv\in V_\lambda$, there exists a compact ball $B\subset V$ such that
for all $i\in\NN$,
\[
B\supset a_i\Omega\cdot\bfv=(a_i\Omega a_i\inv)a_i\cdot\bfv=
\lambda(a_i)(a_i\Omega a_i\inv)\bfv.
\]
Since $\lambda(a_i)\to c$ for some $c>0$ and $\cup_{i\in\NN} a_i\Omega
a_i\inv=U^+$, we have $U^+\cdot \bfv\subset c\inv B$. Since $U^+$ acts
on $V$ by unipotent linear transformations, we obtain that
$U^+\cdot\bfv=\bfv$. Thus $U^+$ acts trivially on $V$. Since the
kernel of $G$ action on $V$ is a normal subgroup of $G$ containing
$U^+$, it is equal to $G$ by our assumption. This contradicts our
hypothesis in the lemma that the action of $G$ is nontrivial. This
proves that $W\subset\sum_{\lambda\in\Lambda_+} V_\lambda$, and the
conclusion of the lemma follows.  \qed

\begin{corollary}  \label{cor:discrete:infinity}
Consider a continuous representation of $G$ on a finite dimensional
vector space $V$ with a Euclidean norm. Let $L=\{\bfv\in
V:G\cdot\bfv=\bfv\}$. Let $\{\bfv_i\}$ be a discrete subset of $V$
contained in $V\setminus L$. Then for any nonempty open set
$\Omega\subset U^+$,
\begin{equation} \label{eq:far}
\sup_{\omega\in\Omega}\|a_i\omega\cdot\bfv_i\|\to\infty \qquad
\mbox{as $i\to\infty$.}
\end{equation}
\end{corollary}

\proof Let $L^\prime$ be the sum of all $G$-invariant irreducible
subspaces of $\dim\geq 2$.  After passing to a subsequence, one of the
following holds:
\[
\mbox{(A) }~~\|\pr_L(\bfv_i)\|\to\infty, \qquad \mbox{or} \qquad
\mbox{(B) }~~\inf_{i\in\NN}\|\pr_{L^\prime}(\bfv_i)\|>0.
\]
If (A) holds then equation~\ref{eq:far} is obvious. If (B) holds, then
there exists an irreducible $G$-subspace $V_1\subset L^\prime$ such
that $\inf_{i\in\NN}\|\pr_{V_1}(\bfv_i)\|>0$. Therefore, without loss
of generality, by replacing $\{\bfv_i\}$ by $\{\pr_{V_1}(\bfv_i)\}$
and $V$ by $V_1$ we may assume that $G$ acts nontrivially and
irreducibly on $V$ and $\inf_{i\in\NN}\|\bfv_i\|>0$.

Let $\omega_0\in\Omega$. Then
$\inf_{i\in\NN}\|\omega_0\bfv_i\|>0$. Therefore replacing $\{\bfv_i\}$
by $\{\omega_0\bfv_i\}$ and $\Omega$ by $\Omega\omega_0\inv$, we may
assume that $e\in\Omega$.

Let $W=\{\bfv\in V:U^+\cdot\bfv=\bfv\}$. By
lemma~\ref{lemma:nilpotent}, there exists $c>0$ such that for all
$i\in\NN$,
\[
\sup_{\omega\in\Omega}\|\pr_W(\omega\cdot\bfv_i)\| 
\geq c\|\bfv_i\|\geq c\cdot\inf_{j\in\NN}\|\bfv_j\|.
\]
Since $\inf_{j\in\NN}\|\bfv_j\|>0$, by lemma~\ref{lemma:take:far},
\[
\sup_{\omega\in\Omega}\|a_i\cdot\omega\bfv_i\|\geq
\sup_{\omega\in\Omega}\|a_i\cdot
\pr_W(\omega\cdot\bfv_i)\|\to\infty\qquad \mbox{as $i\to\infty$}.
\]
\qed


\section{Proofs of the main results}

\subsection*{Translates of horospherical patches} 

\subsection*{\it Proof of theorem~\ref{thm:translate:horosphere}.}

Since $U^+$ is $\sigma$-compact, without loss of
generality we may assume that $\supp(\lambda)$ is compact.  Let
$\la{U}^+$ denote the Lie algebra of $U^+$. We identify $\la{U}^+$
with $\RR^m$ ($m=\dim\la{U}^+$). Let $B$ be a ball in $\la{U}^+$
around the origin such that $\supp(\lambda)\subset\exp(B)$. Let $\nu$
be the restriction of the Lebesgue measure on $B$. By our hypothesis,
$\lambda$ is absolutely continuous with respect to $\exp_\ast(\nu)$,
denoted by $\lambda\abscont\exp_\ast(\nu)$. 

For each $i\in\NN$, define $\Theta_i:\RR^m\to G\subset L$ as
$\Theta_i(\bft)=a_i\exp(\bft)$, $\forall \bft\in
\RR^m\cong\la{U}^+$. Since $\la{U}^+$ is a nilpotent Lie algebra,
there exists $d\in\NN$ such that $\Theta_i\in\cP_{d,m}(L)$, $\forall
i\in\NN$.

\begin{claim} \label{claim:return}
Given $\delta>0$ there exists a compact set $K\subset L/\Lambda$ such
that  
\[
(a_i\pi_\ast(\lambda))(K)>1-\delta, \qquad \forall i\in\NN. 
\]
\end{claim}

Suppose that the claim fails to hold. Since
$\lambda\abscont\exp_\ast(\nu)$, there
exists an $\epsilon>0$ such that for any compact set $K\subset
L/\Lambda$, 
\[
\frac{1}{\nu(B)}(\Theta_i)_\ast(\nu)(K)<1-\epsilon, \qquad \mbox{
for $i$ in a subsequence}.
\]
We apply theorems~\ref{thm:compact:compact} and
\ref{thm:return:strong} for the Lie group $L$, the lattice $\Lambda$,
and the polynomial maps $\Theta_i\in\cP_{d,m}(L)$, $\forall
i\in\NN$. Then by passing to a subsequence, there exists a continuous
representation of $L$ on a finite dimensional vector space $V$ with a
Euclidean norm and a nonzero vector $\bfp\in V$ such that the
following holds: (1) the orbit $\Gamma\cdot\bfp$ is discrete (see
proposition~\ref{prop:discrete}), and (2) for each $i\in\NN$ there
exists $\bfv_i\in\Gamma\cdot\bfp$ such that
\begin{equation}  \label{eq:contracts}
\sup_{\omega\in\exp(B)}\|a_i\omega\cdot\bfv_i\|\to 0\qquad 
\mbox{ as $i\to\infty$}.
\end{equation}

After passing to a subsequence, we may assume that
$G\cdot\bfv_i\neq\bfv_i$, $\forall i\in\NN$. Then
corollary~\ref{cor:discrete:infinity} contradicts
equation~\ref{eq:contracts}. This proves the claim.

By claim~\ref{claim:return}, after passing to a subsequence, we may
assume that the sequence $a_i\cdot\pi_\ast(\lambda)\to\mu$ as
$i\to\infty$, where $\mu$ is a probability measure on $L/\Lambda$.

\begin{claim} \label{claim:invariant}
The measure $\mu$ is $U^+$-invariant.
\end{claim}

To prove the claim, let $u\in U^+$. Then for all $i\in\NN$,
\begin{equation} \label{eq:ua=aui}
u(a_i\pi_\ast(\lambda))=a_i(u_i\pi_\ast(\lambda))=a_i\pi_\ast(u_i\lambda),
\end{equation}
where $u_i=a_i\inv u a_i\in U^+$. Note that $u_i\to e$ as
$i\to\infty$.

Let $\eta$ be a Haar measure on $U^+$. Since $\lambda\abscont\eta$,
there exists $h\in L^1(U,\eta)$ such that $d\lambda=h\,d\eta$.  
Now for any bounded continuous function $f$ on $L/\Lambda$,
\begin{equation}  \label{eq:zero}
\begin{array}{lcl}
&\ & 
\left|\int f \,d[a_i\pi_\ast(u_i\lambda))] - 
\int f\,d[a_i\pi_\ast(\lambda)]\right| 
\\ &=& 
\left|\int_{U^+} f(a_i\pi(u_i\omega))\,d\lambda(\omega) -
\int_{U^+}f(a_i\pi(\omega))\,d\lambda(\omega)\right| 
\\ &=&
\left|\int_{U^+}f(a_i\pi(u_i\omega)h(\omega)\,d\eta(\omega) -
\int_{U^+}f(a_i\pi(\omega))h(\omega)\,d\eta(\omega)\right|
\\&=&
\left|\int_{U^+} f(a_i\pi(\omega))h(u_i\inv\omega)\,d\eta(\omega) -
\int_{U^+}f(a_i\pi(\omega))h(\omega)\,d\eta(\omega)\right|
\\&\leq&
\sup|f|\cdot
\int_{U^+}|h(u_i\inv\omega)-h(\omega)|\,d\eta(\omega) 
\to 0 \qquad \mbox{as $i\to\infty$,}
\end{array}
\end{equation}
because the left regular representation of of $U^+$ on $L^1(U^+,\eta)$
is continuous.

Since $a_i\pi_\ast(\lambda)\to\mu$ as $i\to\infty$, by
equation~\ref{eq:zero}, we get $a_i\pi_\ast(u_i\lambda))\to\mu$ as
$i\to\infty$. Therefore by equation~\ref{eq:ua=aui}, $u\mu=\mu$. This
completes the proof of the claim.

In view of claim~\ref{claim:invariant}, we apply
theorem~\ref{thm:mu_H} to $W=U^+$.  Then there exists a closed
subgroup $H$ of $L$ in the collection $\cH_\Lambda$, such that
\[
\mu(\pi(S_L(H,U^+))=0 \mbox{\qquad and\qquad }\mu(\pi(N_L(H,U^+)))>0.
\]
Let a compact set $C\subset\pi(N_L(H,U^+))\setminus\pi(S_L(H,U^+))$ be
such that $\mu(C)>0$. Since $\lambda\abscont\exp_\ast(\nu)$, there
exists $\epsilon>0$ such that for any Borel measurable set $E\subset
U^+$,
\begin{equation}  \label{eq:lambda<<nu}
\frac{1}{\nu(B)}\exp_\ast(\nu)(E)<\epsilon \Rightarrow
\lambda(E)<\mu(C)/2.
\end{equation}

Let the finite dimensional vector space $V_L$ and the unit vector
$\bfp_H\in V_L$ be as described in notation~\ref{not:V}, for $L$ in
place of $G$ there. We apply theorem~\ref{thm:general-avoid} for
$\epsilon>0$, $d\in\NN$, and $m\in\NN$ chosen as above, and the
compact set $C\subset\pi(N_L(H,U^+))\setminus \pi(S_L(H,U^+))$ as
above. Then there exists a relatively compact set $\Phi\subset V_L$
and an open neighbourhood $\Psi$ of $C$ in $L/\Lambda$ such that for
each $i\in\NN$, applying the theorem to $\Theta_i$ in place of
$\Theta$, one of the following conditions holds:
\begin{enumerate}
\item
There exists $\bfv_i\in\Lambda\cdot\bfp_H$ such that 
\[
a_i\exp(B)\cdot\bfv_i\subset\Phi.
\]
\item
\[
\frac{1}{\nu(B)}\nu(\{\bft\in B:\pi(a_i\exp(\bft))\in\Psi\})<\epsilon.
\]
\end{enumerate}

Since $a_i\pi_\ast(\lambda)\to\mu$ as $i\to\infty$ and $\Psi$ is a
neighbourhood of $C$, there exists $i_0\in\NN$ such that
$\lambda(\pi\inv(a_i\inv\Psi)\cap U^+)>\mu(C)/2$ for all $i\geq i_0$.
Therefore by equation~\ref{eq:lambda<<nu}, condition~1 must hold for
all $i\geq i_0$. Now by passing to a subsequence, there exists
$\bfv_i\in\Lambda\cdot\bfp_H$ for each $i\in\NN$ such that
\begin{equation}  \label{eq:bounded}
a_i\exp(B)\cdot\bfv_i\subset\Phi.
\end{equation}
By proposition~\ref{prop:discrete}, the sequence $\{\bfv_i\}$ is
discrete. By corollary~\ref{cor:discrete:infinity} and
equation~\ref{eq:bounded}, there exists $i_0\in\NN$ such that 
$G\cdot\bfv_{i_0}=\bfv_{i_0}$. Let $\gamma\in\Lambda$ such that
$\bfv_{i_0}=\gamma\bfp_H$. Then
\[
G\cdot\gamma\cdot\bfp_H=\gamma\cdot\bfp_H.
\]
Thus $G\subset \gamma N_L^1(H)\gamma\inv$. But $\pi(N_L^1(H))$ is
closed in $L/\Lambda$ by proposition~\ref{prop:discrete}, and $\pi(G)$
is dense in $L/\Lambda$. Therefore we conclude that $H$ is a normal
subgroup of $L$. Since $N_L(H,U^+)\supset C\neq\emptyset$, this
implies in particular that $U^+$ is contained in $H$. Thus $U^+\subset
G\cap H$ and $G\cap H$ is normal in $G$. Therefore by our hypothesis
$G\cap H=G$, or in other words $G\subset H$. Again since $\pi(G)$ is
dense in $L/\Lambda$, we have $H=L$. Therefore
$\mu(\pi(S(L,U^+)))=0$. Hence by theorem~\ref{thm:mu_H}, we have that
$\mu$ is $L$-invariant. This completes the proof of the theorem.  \qed

\subsection*{Translates of orbits of symmetric subgroups}

First we make some observations. For the results stated below, let
$(U,\nu_1)$ and $(V,\nu_2)$ be locally compact second countable spaces
with Borel measures.

\begin{proposition} \label{prop:decomposition}
Let $\lambda$ be a Borel probability measure on $U\times V$ which is
absolutely continuous with respect to $\nu_1\times\nu_2$, denoted by
$\lambda\abscont\nu_1\times\nu_2$. Then there exists a probability
measure $\lambda_1\abscont\nu_1$ on $U$, and for almost all $u\in
(U,\lambda_1)$, there exists a probability measure
$\lambda_u\abscont\delta_u\times\nu_2$ on $\{u\}\times V$, where
$\delta_u$ is the point mass at $\{u\}$, such that the following
holds: For any bounded continuous function $f$ on $U\times V$, the map
$u\mapsto{\scriptsize \int_{\{u\}\times V} f\, d\lambda_u}$ is
$\lambda_1$-measurable, and
\[
\int_{U\times V} f\,d\lambda=\int_U \left(\int_{\{u\}\times V} f\,
d\lambda_u\right)d\lambda_1(u).
\] 
\end{proposition}

\proof Let $h=d\lambda/d(\nu_1\times\nu_2)\geq 0$ be the Radon-Nikodym
derivative. For any $u\in U$, put $\alpha(u)=\int_V
h(u,v)\,d\nu_2(v)$. Let $C=\{u\in U:\alpha(u)>0\}$. Let $\lambda_1$ be
the restriction of $\nu_1$ to $C$. For almost any $u\in (U,\lambda_1)$,
let $\lambda_u$ be the Borel measure on $\{u\}\times V$ such
that $d\lambda_u/d[\delta_u\times\nu_2]=h(u,\cdot)/\alpha(u)$. Now the
conclusion of the proposition follows from Fubini's theorem.  
\qed

For the propositions stated below, let $G$ be a locally
compact topological group acting continuously on a locally compact
space $X$. Let $\{a_i\}$ be a sequence in $G$ and $\mu$ a Borel
probability measure on $X$.

\begin{proposition} \label{prop:invariant}
Let $\lambda$ be a probability measure on $X$ such that
$a_i\lambda\to\mu$  as $i\to\infty$. Let $b\in G$ such that
$\cl{\{a_iba_i\inv:i\in\NN\}}$ is compact. If $\mu$
is $G$-invariant, then $a_i(b\lambda)\to\mu$  as $i\to\infty$.
\end{proposition}

\proof
First observe that there is no loss of generality in passing to a
subsequence. Therefore we may assume that
$a_iba_i\inv\to g$ for some $g\in G$. Now 
\[
a_i(b\lambda)=(a_iba_i\inv)(a_i\lambda)\to g\mu \qquad \mbox{as
$i\to\infty$}.
\]
Since $g\mu=\mu$, the proof is complete. 
\qed

For the next two propositions, assume that $G$ contains the spaces $U$
and $V$. Fix $x_0\in X$, and let $\rho:U\times V\to X$ be the map
given by $\rho(u,v)=uvx_0$, $\forall (u,v)\in U\times V$.

\begin{proposition} \label{prop:contract}
Let the notation be as in
proposition~\ref{prop:decomposition}. Suppose that for almost all
$u\in (U,\lambda_1)$, we have $a_i\rho_\ast(\lambda_u)\to\mu$  as
$i\to\infty$. Then $a_i\rho_\ast(\lambda)\to\mu$  as
$i\to\infty$.
\end{proposition}

\proof
Let $f$ be bounded continuous function on $X$. Then 
\[
\begin{array}{lcl}
\int_X f \,d[a_i\rho_\ast(\lambda)]
&=& \int_{U\times V} f(a_i\rho(\omega))\,d\lambda(\omega) 
\\ &=& 
\int_{V} d\lambda_1(u)\cdot
\int_{\{u\}\times V} f(a_i\rho(\omega))\,d\lambda_u(\omega) 
\\ &=&
\int_{V} d\lambda_1(u) \cdot \int_X f\,d[a_i\rho_\ast(\lambda_u)]
\\ &\to& 
\int_{V}\, d\lambda_1(u) \cdot \int_X f\,d\mu \qquad 
\mbox{as $i\to\infty$}
\\ &=& \int_X f\, d\mu. 
\end{array}
\]
\qed

By similar arguments we obtain the following result.

\begin{proposition}  \label{prop:same}
Suppose that $a_i u a_i\inv\to e$ as $i\to\infty$ for all $u\in
U$. Then as $i\to\infty$, 
\[
a_i\rho_\ast(\nu_2)\to\mu \Leftrightarrow
a_i\rho_\ast(\nu_1\times\nu_2)\to\mu. 
\]
\end{proposition} 

\subsection*{\it Proof of corollary~\ref{cor:symmetric}.~} 

Using the results in \cite[Section 7.1]{Sch:book} there exist an
$\RR$-split torus $A\subset G$ and a maximal compact subgroup $K$ of
$G$ such that the following holds: (1)~$\sigma(a)=a\inv$, $\forall
a\in A$, (2)~the set of real roots of $A$ for the adjoint action on
the Lie algebra of $G$ forms a root system, and (3)~$G$ admits a
decomposition $G=K\cl{A}^+ S$, were $\cl{A}^+$ denotes the closure of
the positive Weyl chamber with respect to a system $\Delta$ of simple
roots on $A$.

Using this decomposition and by passing to a subsequence, without loss
of generality we may assume the following: (1) $g_i=a_i\in\cl{A}^+$
for all $i\in\NN$; (2) $\{a_i\}_{i\in\NN}$ has no convergent
subsequence, (because otherwise $G_1=\{e\}$ and $\pi(e)$ cannot be
dense in $L/\Lambda$); and (3) for any $\alpha\in\Delta$, either
$\sup_{i\in\NN}\alpha(a_i)<\infty$ or $\alpha(a_i)\to\infty$ as
$i\to\infty$. 

For the rest of the proof, consider the notation~\ref{not:normal}.

Let $G_1$ be the smallest closed normal subgroup of $G$ containing
$U^+$.  Then it is straightforward to verify that the projection of
$\{a_i\}$ on $G/G_1$ is relatively compact. Therefore by our
hypothesis, $\cl{\pi(G_1)}=L/\Lambda$. 

Take any $g_0\in S$ and define $\rho(h)=\pi(hg_0)$ for all $h\in L$.
Since any closed connected normal subgroup of $G_1$ is also normal in
$G$, we can apply theorem~\ref{thm:translate:horosphere} to $G_1$ in
place of $G$ and $\rho$ in place of $\pi$. Then for any probability
measure $\nu$ on $U^+$ which is absolutely continuous with respect to
a Haar measure on $U^+$, we have
\begin{equation}  \label{eq:bi.lambda}
a_i\rho_\ast(\nu)\to \mu_L, \qquad \mbox{as $i\to\infty$}.
\end{equation}

Since $\sigma(a)=a\inv$ ($\forall a\in A$), for any $X\in\la{U}^+$, we
have $\sigma(X)\in\la{U}^-$ and $X+\sigma(X)\in \la{S}$. Also
$\sigma(\la{Z})=\la{Z}$. Now by equation~\ref{eq:G=U-ZU+},
\begin{equation}  \label{eq:U-S=U0U+}
\la{U}^-\oplus\la{S}=\la{U}^-\oplus(\la{S}\cap\la{Z})\oplus\la{U}^+.
\end{equation}
Then by implicit function theorem, there exist relatively compact
neighbourhoods $\Omega^-$, $\Omega^0$, $\Omega^+$ and $\Phi$ of $e$ in
$U^-$, $(Z\cap S)U^-$, $U^+$ and $S$, respectively, such that for any
open set $\Psi$ of $\Phi$, we have that $\Omega^-\Psi$ is an
open subset of $\Omega^0\Omega^+$. Also we may assume that under the
multiplication map $\Omega^-\times\Phi\cong \Omega^-\Phi$ and
$\Omega^0\times\Omega^+\cong\Omega^0\Omega^+$. 

Let $\nu_-$ and $\nu'$ be probability measures obtained by restricting
Haar measures of $U^-$ and $S$ to $\Omega^-$ and $\Psi$,
respectively. Then $\lambda=\nu_-\times\nu'$ is a smooth measure on
$\Omega^-\times\Psi$. By choosing $\Psi$ small enough, we can ensure
that $\rho_\ast(\nu')$ is a multiple of $\mu_S$ restricted to
$\rho(\Psi)$. Since $g_0\in S$ chosen in the definition of $\rho$ is
arbitrary and since there is enough flexibility in the choices of
$\Phi$ and $\Psi$, to prove that $a_i\mu_S\to\mu_L$, it is enough to
show that $a_i\rho_\ast(\nu')\to\mu_L$ as $i\to\infty$.

By proposition~\ref{prop:same}, as $i\to\infty$,
$a_i\rho_\ast(\nu')\to\mu_L$ if and only if
$a_i\rho_\ast(\lambda)\to\mu_L$. Therefore to complete
the proof of the corollary, it is enough to show the following.

\begin{claim}  \label{claim:bi.lambda}
As $i\to\infty$, $a_i\rho_\ast(\lambda)\to\mu_L$. 
\end{claim}

Since $\Omega^-\Psi\subset\Omega^0\Omega^+$, $\lambda$ can be treated
as a measure on $\Omega^0\times\Omega^+$. Let $\nu_1$ and $\nu_2$ be
the probability measures obtained by restricting the Haar measures on
$(Z\cap S)U^-$ and $U^+$ to $\Omega^0$ and $\Omega^+$, respectively.
Since $\lambda$ is a smooth measure, $\lambda\abscont\nu_1\times\nu_2$
(see equation~\ref{eq:U-S=U0U+}). Decompose $\lambda$ as in
proposition~\ref{prop:decomposition}. Then for almost all
$\omega\in(\Omega^0,\lambda_1)$, we have
$\lambda_\omega\abscont\omega\nu_2$. Put
$\nu_\omega=\omega\inv\lambda_\omega$. Then $\nu_\omega\abscont\nu_2$.
Hence by equation~\ref{eq:relcpt}, equation~\ref{eq:bi.lambda},
and proposition~\ref{prop:invariant},
\[
a_i\rho_\ast(\lambda_\omega)=a_i(\omega\rho_\ast(\nu_\omega))\to\mu_L
\qquad \mbox{as $i\to\infty$}.
\]
Now by proposition~\ref{prop:contract}, $a_i\pi_\ast(\lambda)\to\mu_L$
as $i\to\infty$. This completes the proof of the claim, and also the
proof of the corollary.  \qed

\subsection*{Continuous $G$-equivariant factors of $G/P\times L/\Lambda$} 

First we recall the following result from
\cite[Section~2]{Dani:factors}.

\begin{proposition}[Dani]  \label{prop:length:parabolic}
Let $G$ be a semisimple group with finite center and
$\RR$-$\rank(G)\geq2$. Let $P$ be a parabolic subgroup of $G$. Then
given $g\in G\setminus P$, there exist $k\in\NN$ ($k<\RR$-$\rank(G)$),
elements $g_1,\ldots,g_{k+1}$ in $G$, and one-parameter unipotent
subgroups $\{u_1(t)\},\ldots,\{u_k(t)\}$ of $G$ contained in $P$ such
that the following holds:
\begin{enumerate}
\item
$g_1=g$, $g_k\not\in P$, and $g_{k+1}=e$.
\item 
For each $i=1,\ldots,k$, 
\[
u_i(t)g_iP\to g_{i+1}P \mbox{ in $G/P$ as $t\to\infty$}.
\]
\item
There exists a semisimple element $a$ of $G$ in $g_kPg_k\inv\cap P$ such
that if $U^+$ is the associated horospherical subgroup then  
$U^+\subset g_kPg_k\inv\cap P$, and if $G_1$ denotes the smallest normal
subgroup of $G$ containing $U^+$, then $\RR$-$\rank(G/G_1)\leq 1$.
\end{enumerate}
\end{proposition}

\proof
Apply \cite[Corollary~2.3]{Dani:factors} iteratively. Also use the
proofs of \cite[Corollary~2.6 and Lemma~2.7]{Dani:factors}.
\qed

Now we obtain the analogue of \cite[Lemma~1.4]{Dani:factors} by using
theorem~\ref{thm:parabolic} in place of
\cite[Lemma~1.1]{Dani:factors}.  Also we use the recurrence conclusion
of theorem~\ref{thm:Ratner} of Ratner in place of
\cite[Lemma~1.6]{Dani:factors}.

\begin{proposition}  \label{prop:factors}
Let the notation and assumptions be as in theorem~\ref{thm:G/P}.  Let
$x,y\in L/\Lambda$ and $g\in G\setminus P$. If $\phi(x,gP)=\phi(y,P)$,
then there exists a parabolic subgroup $Q$ containing $\{g\}\cup P$
such that $\phi(z,P)=\phi(z,qP)$ for all $z\in\cl{Gx}$ and $q\in Q$.
Moreover, $\phi(y,P)=\phi(y,qP)$ for all $q\in Q$. 
\end{proposition}

\proof
Let $k\in\NN$, elements $g_1,\ldots,g_{k+1}$ in
$G$, the one-parameter unipotent subgroups $\{u_i(t)\}$
contained in $P$, and a semisimple element $a$ of $G$ and the
associated expanding horospherical subgroup $U^+$ be as in
proposition~\ref{prop:length:parabolic}. For each $i=1,\ldots, k$,
by Ratner's theorem~\ref{thm:Ratner} applied to the diagonal action of
$\{u_i(t)\}$ on $L/\Lambda\times L/\Lambda$, there exists a sequence
$t_n\to\infty$ such that $(u_i(t_n)x,u_i(t_n)y)\to (x,y)$ as $n\to\infty$.
Now for any $i\in \{1,\ldots,k\}$,
\[
\phi(x,g_iP)=\phi(y,P) \Rightarrow
\phi(u_i(t_n)x,u_i(t_n)g_iP)=\phi(u_i(t_n)y,P), \forall n\in\NN. 
\]
In the limit as $n\to\infty$,  we get
$\phi(x,g_{i+1}P)=\phi(y,P)$. Since $g_1=g$, by  induction on $i$, we get
that $\phi(x,g_iP)=\phi(y,P)$ for all $1\leq i\leq k+1$. 

In particular, since $g_{k+1}=e$,
\[
\phi(x,g_k P)=\phi(y, P)=\phi(x,P).
\]
Since $F=\{a^n:n\in\NN\}\cdot U^+\subset g_k P g_k\inv\cap P$, we have
that
\[
\phi(fx,g_kP)=\phi(fx,P),\,\forall f\in F.
\] 

Let $G_1$ be the smallest closed normal subgroup of $G$ containing
$U^+$.  Then by the choice of $a$ as in
Proposition~\ref{prop:length:parabolic}, $\RR$-$\rank(G/G_1)\leq
1$. Therefore by the hypothesis in theorem~\ref{thm:G/P},
$\cl{G_1x}=\cl{Gx}$. By theorem~\ref{thm:Ratner}, $\cl{Gx}$ is an
orbit of a closed subgroup, say $L'$, of $L$ containing $G$, and the
stabilizer of $x$ in $L'$, say $\Lambda'$, is a lattice in
$L'$. Applying theorem~\ref{thm:parabolic} to $L'$ and $\Lambda'$ in
place of $L$ and $\Lambda$, respectively, we conclude that
$\cl{Fx}=\cl{Gx}$.  Thus
\[
\phi(z,g_kP)=\phi(z,P),\qquad \forall z\in\cl{G_1x}=\cl{Gx}.
\]

Put 
\begin{equation}  \label{eq:Q}
Q=\{h\in G:\phi(z,fhP)=\phi(z,fP),\ \forall z\in\cl{Gx}\mbox{ and }\forall
f\in G\}. 
\end{equation}
Then $Q$ is a closed subgroup of $G$ containing $P\cup\{g_k\}$. Since
$g_k\not\in P$, 
\begin{equation}  \label{eq:dim:increase}
Q\neq P.
\end{equation}

Now if $g\not\in Q$, then replacing $P$ by $Q$ and $L/\Lambda$ by
$\cl{Gx}$, we repeat the whole argument. Note that by definition the
new set given by equation~\ref{eq:Q} still turns out to be same as
$Q$.  This fact contradicts the new
equation~\ref{eq:dim:increase}. This completes the proof. 
\qed

\subsection*{\it Proof of theorem~\ref{thm:G/P}.}

Define the equivalence relation
\[
R=\{(x,y)\in L/\Lambda\times L/\Lambda:\phi(x,gP)=\phi(y,gP)\mbox{ for
some $g\in G$}\}
\]
on $L/\Lambda$. Clearly $R$ is a closed subset of $L/\Lambda\times
L/\Lambda$ invariant under the diagonal action of $G$. Let $X$ be the
space of equivalence classes of $R$ and let $\phi_1:L/\Lambda\to X$ be
the map taking any element of $L/\Lambda$ to its equivalence
class. Equip $X$ with the quotient topology. Then $X$ is a locally
compact Hausdorff space.

For any $x\in L/\Lambda$, put
\[
\cQ(x)=\{h\in G:\phi(x,gP)=\phi(x,ghP),\ \forall g\in G\}.
\]
Observe that $\cQ(x)$ is a closed subgroup of $G$ containing $P$ and
for any $y\in\cl{Gx}$, we have $\cQ(y)\supset\cQ(x)$. Let $x_0\in
L/\Lambda$ such that $\cl{Gx_0}=L/\Lambda$ and put $Q=\cQ(x_0)$. Then
$\cQ(y)\supset Q$ for all $y\in L/\Lambda$.  Since $Q$ is a parabolic
subgroup of $G$, there are only finitely many closed subgroups of $G$
containing $Q$. Therefore the set $X_Q:=\{x\in L/\Lambda:\cQ(x)=Q\}$
is open in $L/\Lambda$. Also $X_Q$ is nonempty and $G$-invariant. Now
since $G$ acts ergodically on $L/\Lambda$, the set $L/\Lambda
\setminus X_Q$ is closed and nowhere dense.

Note that for any $x,y\in L/\Lambda$, if $\phi_1(x)=\phi_1(y)$ then by
proposition~\ref{prop:factors}, we have that $\cQ(x)=\cQ(y)$.  Let
$\rho:L/\Lambda\times G/P\to X\times G/Q$ be the ($G$-equivariant) map
defined by $\rho(x,gP)=(\phi_1(x),gQ)$ for all $x\in L/\Lambda$ and
$g\in G$. Then there exists a uniquely defined map $\psi:X\times
G/Q\to Y$ such that $\phi=\psi\circ\rho$. It is straightforward to
verify that $\psi$ is continuous and $G$-equivariant.

Take any $x\in X_Q$, $y\in L/\Lambda$, and $g,h\in G$ such that
$\phi(x,ghP)=\phi(y,gP)$. Then $\phi_1(y)=\phi_1(x)$, and hence
$h\in\cQ(y)=\cQ(x)=Q$. This proves that $\psi$ restricted to
$\phi_1(X_Q)\times G/Q$ is injective and $y\in X_Q$.  

Now if $Y$ is a locally compact second countable space and $\phi$ is
surjective, then using Baire's category theorem for Hausdorff locally
compact second countable spaces, one can show that $\phi$ is an open
map. This completes the proof of the theorem.  
\qed

\subsection*{Continuous $G$-equivariant factors of $L/\Lambda$}

\subsection*{\it Proof of theorem~\ref{thm:factors}.~} 

Define $\Lambda_1=\{h\in L:\phi(gh\Lambda)=\phi(g\Lambda),\,\forall
g\in L\}$. Then $\Lambda_1$ is a closed subgroup of $L$ containing
$\Lambda$.  Since $G$-acts ergodically on $L/\Lambda$, $\Ad(\Lambda)$
is Zariski dense in $\Ad(L)$
(see~\cite[Theorem~2.3]{Sh:equidistribution}). Therefore $\Lambda_1^0$
is a normal subgroup of $L$. Let $\Lambda_1^\prime$ be the largest
subgroup of $\Lambda_1$ which is normal in $L$. In view of
\ref{rem:affine}~(3), replacing $L$ by $L/\Lambda_1^\prime$, $\Lambda$
by $\Lambda_1/\Lambda_1^\prime$, and $G$ by its image in
$L/\Lambda_1'$, without loss of generality we may assume that
$\Lambda'_1=\{e\}$ and $\Lambda_1=\Lambda$.

Define the equivalence relation
\[
R=\{(x,y)\in L/\Lambda\times L/\Lambda:\phi(x)=\phi(y)\}
\]
on $L/\Lambda$. Then $R$ is closed and $\Delta(G)$-invariant, where
$\Delta:L\to L\times L$ denotes the diagonal embedding of $L$ in $L\times
L$. 

Let
\[
 K=\{\tau\in\Aff(L/\Lambda):(z,\tau(z))\in R\mbox{ and }
\tau(gz)=g\tau(z),\,\forall z\in L/\Lambda,\,\forall g\in G\}
\]
and
\[
X_1=\{x\in L/\Lambda:\cl{Gx}=L/\Lambda\}.
\]
Note that  $X_1\neq\emptyset$, since $G$ acts ergodically on $L/\Lambda$. 

\begin{claim} \label{claim:tau}
Let $(x,y)\in R$. If $x\in X_1$, then $y\in X_1$ and there exists
$\tau\in K$ such that $y=\tau(x)$.
\end{claim}

The claim is proved as follows. Since $\Delta(G)$ is generated by
one-parameter unipotent subgroups of $L\times L$, by Ratner's
theorem~\ref{thm:Ratner} there exists a closed subgroup $F$ of
$L\times L$ containing $\Delta(G)$ such that
\[
\cl{\Delta(G)\cdot(x,y)}=F\cdot(x,y)
\]
and $F\cdot(x,y)$ admits an $F$-invariant probability measure, say
$\lambda$. 

Let $p_i:L\times L\to L$ denote the projection on the $i$-th
coordinate, where $i=1,2$. Then $(\pi\circ p_1)_\ast(\lambda)$ is a
$p_1(F)$-invariant probability measure on $p_1(F)x$. Hence the orbit
$p_1(F)x$ is closed (see~\cite[theorem~1.13]{Rag:book}).  Since
$G\subset p_1(F)$ and $\cl{Gx}=L/\Lambda$, we have that
$p_1(F)=L$. Let $N_1=p_1(F\cap\ker(p_2))$. Then $N_1$ is a normal
subgroup of $p_1(F)=L$ and $(N_1z,w)\subset R$ for all $(z,w)\in
F\cdot(x,y)$.  Therefore $N_1\subset\Lambda_1'=\{e\}$. Thus
$F\cap\ker(p_2)=N_1\times\{e\}=\{e\}$. Now since $p_1(F)=L$ and
$p_2|_F$ is injective, $\dim(p_2(F))=\dim(L)$. Since $L$ is
connected, $p_2(F)=L$. Thus $p_2|_F$ is an isomorphism.

Now $\cl{Gy}\supset p_2(F)y=L/\Lambda$. Hence $y\in X_1$. Now
interchanging the roles of $x$ and $y$ in the above argument, we
conclude that $p_1|_F$ is an isomorphism. Let
$\sigma=p_2\circ(p_1|_F)\inv$. Then $\sigma\in\Aut(L)$ and
\[
F=\{(g,\sigma(g))\in L\times L:g\in L\}.
\]

Thus $(gx,\sigma(g)y)\in R$ for all $g\in L$. Now for any $\delta\in
L$, if $\delta x=x$, then $(gx,\sigma(g)\sigma(\delta)y)\in R$ for all
$g\in L$. Let $h\in L$ such that $y=h\Lambda$. Then
$\phi(\sigma(g)h\Lambda)=\phi(\sigma(g)\sigma(\delta)h\Lambda)$ for
all $g\in L$. Since $\sigma(L)h=L$, we conclude that
$h\inv\sigma(\delta)h\in\Lambda_1$. Now since $\Lambda_1=\Lambda$, we
have that $\sigma(\delta)y=y$. Therefore the map $\tau:L/\Lambda\to
L/\Lambda$, given by $\tau(gx)=\sigma(g)y$ for all $g\in L$, is well
defined and $\tau\in\Aff(L/\Lambda)$.

Therefore
\[
F\cdot(x,y)=\{(z,\tau(z)):z\in L/\Lambda\}. 
\]
Since $\Delta(G)\subset F$, we have that $\sigma(g)=g$ 
and hence $\tau(gz)=gz$, for all $g\in G$. Thus $\tau\in K$, and the
proof of the claim is complete. 

\begin{claim} \label{claim:cK:compact}
The group $ K$ is compact.
\end{claim}

We prove the claim as follows. Clearly, $K$ is a closed subset of
$\Aff(L/\Lambda)$, and hence it is locally compact. Let $\mu_L$ denote
the $L$-invariant probability measure on $L/\Lambda$. Then
$\mu_L(X_1)=1$. For any $x\in X_1$, if $y\in\cl{K\cdot x}$ then
$(x,y)\in R$, and by claim~\ref{claim:tau} there exists $\tau\in K$
such that $y=\tau(x)$. Thus $ K\cdot x$ is closed for all $x\in
X_1$. Therefore by Hedlund's Lemma and the ergodic decomposition of
$\mu_L$ with respect to the action of $K$ on $L/\Lambda$, we have that
almost all $K$-ergodic components are supported on closed $
K$-orbits. Thus for almost all $x\in L/\Lambda$, the orbit $ K\cdot x$
supports a $K$-invariant probability measure.

For any $x\in L/\Lambda$, put $ K_x=\{\tau\in K:\tau(x)=x\}$. Let
$\xi: K/ K_x\to L/\Lambda$ be the map defined by $\xi(\tau
K_x)=\tau(x)$ for all $\tau\in K$.  Since $\Aff(L/\Lambda)$ acts
continuously on $L/\Lambda$, we have that $\xi$ is a continuous
injective $K$-equivariant map. Let $x\in X_1$ be such that $K\cdot x$
supports a $K$-invariant probability measure. Since $\xi$ is
injective, the measure can be lifted to a $K$-invariant probability
measure on $K/K_x$. Let $\tau\in K_x$. Then for any $g\in G$, we have
$\tau(gx)=g\tau(x)=gx$. Now since $\cl{Gx}=L/\Lambda$, we have that
$\tau(y)=y$ for all $y\in L/\Lambda$. Hence $K_x$ is the trivial
subgroup of $\Aff(L/\Lambda)$. Thus $K$ admits a finite Haar measure.
Hence $K$ is a compact group, and the claim is proved.

Let $\Omega$ be any neighbourhood of $e$ in $Z_L(G)$. Put
\[
R'=\{(x,y)\in R: y\not\in K\cdot\Omega x\}.
\]
Let $X_c$ be the closure of the projection of $R'$ on the first factor
of $L/\Lambda\times L/\Lambda$. Put $X_0=(L/\Lambda)\setminus X_c$.

\begin{claim} 
$X_1\subset X_0$. 
\end{claim}

Suppose the claim does not hold. Then there exists a sequence
$\{(x_i,y_i)\}\subset R'$ converging to $(x,y)\in R$ with $x\in
X_1$. By claim~\ref{claim:tau}, there exists $\tau\in K$ such that
$y=\tau(x)$.  Therefore, after passing to a subsequence, there exists
a sequence $g_i\to e$ in $L$ such that $y_i=\tau(g_ix_i)$ for all
$i\in\NN$.  By the definition of $R'$, $g_i\not\in\Omega\subset
Z_L(G)$ for all $i\in\NN$.  Also $(x_i,g_ix_i)\in R$ for all
$i\in\NN$. By Ratner's theorem~\ref{thm:Ratner}, there exits a
$\Delta(G)$-invariant $\Delta(G)$-ergodic probability measure $\mu_i$
on $(L/\Lambda)\times(L/\Lambda)$ such that
$\cl{\Delta(G)(x_i,g_ix_i)}=\supp(\mu_i)$. Let $h_i\to e$ be a
sequence in $L$ such that $x_i=h_ix$ for all $i\in\NN$.  By
theorem~\ref{thm:limit}, after passing to a subsequence, there exists
a probability measure $\mu$ on $L/\Lambda\times L/\Lambda$ such that
$\mu_i\to\mu$ as $i\to\infty$, and the following holds:
$\supp(\mu)=F\cdot (x,x)$, where $F$ is a closed subgroup of $L\times
L$, and
\begin{equation}  \label{eq:conjugate}
(h_i\inv,h_i\inv g_i\inv)\Delta(G)(h_i,g_ih_i)\subset F,\qquad \forall
i\in\NN.
\end{equation}
In particular, $F\cdot (x,x)\subset R$ and $\Delta(G)\subset F$. Since
$x\in X_1$, we have that $F\supset\Delta(L)$. By an argument as in the
proof of claim~\ref{claim:tau}, we conclude that
$F\cap\ker(p_i)=\{e\}$ for $i=1,2$. Therefore $F=\Delta(L)$. Hence
from equation~\ref{eq:conjugate} we conclude that $g_i\in Z_L(G)$,
which is a contradiction. This completes the proof of the claim, and
the proof of the theorem.  \qed

\bigskip\noindent{\bf Acknowledgements}

\medskip\noindent The author wishes to thank S G Dani, Alex Eskin, G A
Margulis, Shahar Mozes, Marina Ratner and Garret Stuck for fruitful
discussions. The suggestions of the referee have been very valuable in
rectifying some errors and improving the presentation.


\bigskip
\noindent
{\bf Note added in the proof:} Equidistributions of translates of
measures is also considered in a recent preprint (later published as
\cite{EMM:upperbound}) ``Upper bounds and asymptotics in a
quantitative version of the Oppenheim conjecture'' of A~Eskin,
G~A~Margulis and S~Mozes. In the context of the preprint it seems
worthwhile to remark that the method in the present paper can be used
to obtain `uniform versions' of theorem~1.4 and corollary~1.2, as done
in the above mentioned paper, for certain results.

\end{document}